\documentclass{compositio}

\usepackage{amsmath}
\usepackage{dsfont}
\usepackage[all]{xypic}
\xyoption{dvips}
\usepackage{paralist}

\theoremstyle{definition}
\newtheorem{ntn}{Notation}[section]
\newtheorem{dfn}[ntn]{Definition}
\theoremstyle{plain}
\newtheorem{lem}[ntn]{Lemma}
\newtheorem{prp}[ntn]{Proposition}
\newtheorem{thm}[ntn]{Theorem}
\newtheorem{cor}[ntn]{Corollary}
\newtheorem{cnj}[ntn]{Conjecture}
\theoremstyle{remark}

\newtheorem{exa}[ntn]{Example}

\DeclareMathOperator{\gr}{gr}
\DeclareMathOperator{\tr}{tr}
\DeclareMathOperator{\Ann}{Ann}
\DeclareMathOperator{\Der}{Der}
\DeclareMathOperator{\II}{I}
\DeclareMathOperator{\LL}{L}
\DeclareMathOperator{\Sp}{Sp}
\DeclareMathOperator{\TT}{T}
\DeclareMathOperator{\Sing}{Sing} 
\DeclareMathOperator{\VV}{V}

\newcommand{\dd}{\mathrm{d}}
\newcommand{\m}{\mathfrak{m}}
\newcommand{\ol}[1]{{\overline{#1}}}
\newcommand{\ul}[1]{{\underline{#1}}}
\newcommand{\wh}[1]{{\widehat{#1}}}
\newcommand{\xymat}{\SelectTips{cm}{}\xymatrix}
\newcommand{\C}{\mathds{C}}
\newcommand{\Cen}{\mathfrak{C}}
\newcommand{\CC}{\mathrm{C}}

\newcommand{\DD}{\mathfrak{D}}
\newcommand{\D}{\mathcal{D}}
\newcommand{\HH}{\mathrm{H}}

\renewcommand{\O}{\mathcal{O}}
\newcommand{\Q}{\mathds{Q}}
\newcommand{\R}{\mathbf{R}}
\renewcommand{\SS}{\mathfrak{S}}
\newcommand{\U}{\mathcal{U}}
\newcommand{\V}{\mathcal{V}}

\begin{document}

\title[On the formal structure of logarithmic vector fields]
{On the formal structure of logarithmic vector fields}

\author{Michel Granger}
\email{michel.granger@univ-angers.fr}
\address{
Universit\'e d'Angers\\
2 Bd. Lavoisier\\
49045 Angers\\
France}

\author{Mathias Schulze}
\email{mathias.schulze@univ-angers.fr}
\address{
Universit\'e d'Angers\\
2 Bd Lavoisier\\
49045 Angers\\
France}
\thanks{The second author was supported by \'EGIDE.}

\shortauthors{M. Granger and M. Schulze}



\begin{abstract}
In this article, we prove that a free divisor in a three dimensional complex manifold must be Euler homogeneous in a strong sense if the cohomology of its complement is the hypercohomology of its logarithmic differential forms.
F.J.~Calder\'on--Moreno et al.\ \cite{CMNC02} conjectured this implication in all dimensions and 
proved it in dimension two.

We prove a theorem which describes in all dimensions a special minimal system of generators for the module of formal logarithmic vector fields.
This formal structure theorem is closely related to the formal decomposition of a vector field in K.~Saito \cite{Sai71} and is used in the proof of the above result. 

Another consequence of the formal structure theorem is that the truncated Lie algebras of logarithmic vector fields up to dimension three are solvable.
We give an example that this may fail in higher dimensions.
\end{abstract}

\maketitle

\section{Introduction}

Let $X$ be a complex manifold of dimension $n\ge1$.
Let $D\subseteq X$ be a reduced divisor and $\xymat{j:U\ar@{^(->}[r]&X}$ the inclusion of its complement $U:=X\backslash D$ in $X$.
By Grothendieck's comparison theorem \cite[Thm.~2]{Gro66}, the natural morphism
\[
\xymat{
\Omega_X^\bullet(*D)\ar[r]&j_*\Omega_U^\bullet\simeq\R j_*\C_U
}
\]
is a quasi--isomorphism.
This means that the complex $\Omega_X^\bullet(*D)$ of holomorphic
differential forms with meromorphic poles along $D$ can be used to calculate the cohomology of $U$.

Generalizing ideas of P.~Deligne and N.~Katz, K.~Saito \cite{Sai80} defined the subcomplex $\Omega_X^\bullet(\log D)\subseteq\Omega_X^\bullet(*D)$ of holomorphic differential forms with logarithmic poles along $D$.
Unlike $\Omega_X^\bullet(*D)$, $\Omega_X^\bullet(\log D)$ is a complex of coherent $\O_X$--modules.
If $D$ is a normal crossing divisor, $\Omega_X^\bullet(\log D)$ also computes the cohomology of $U$.
This fact plays a crucial role in P.~Deligne's mixed Hodge theory \cite[\S3]{Del71}.
In general, one says that the logarithmic comparison theorem, or simply LCT, holds for $D$ if the inclusion
\[
\xymat{
\Omega_X^\bullet(\log D)\ar@{^(->}[r]&\Omega_X^\bullet(*D)
}
\]
is a quasi--isomorphism.
The characterization of LCT in general is an open problem.

The natural dual of $\Omega_X^1(\log D)$ is the module $\Der_X(-\log D)$ of logarithmic vector fields along $D$.
As recently proposed by K.~Saito, we adopt a notation harmonised with the conventions of algebraic geometry.
For $x\in X$, a vector field $\delta\in\Der_{X,x}$ is contained in $\Der_{X,x}(-\log D)$ if $\delta(f)\in\O_{X,x}\cdot f$ for some, and hence any, (reduced) local equation $f\in\O_{X,x}$ of $(D,x)$.
We shall often use the standard notation $\II(D,x)=\langle f\rangle:=\O_{X,x}\cdot f$ or $(D,x)=\VV(f)$.
K.~Saito \cite{Sai80} introduced the important class of free divisors:
The divisor $D$ is called free if $\Omega_X^1(\log D)$, or equivalently $\Der_X(-\log D)$, is a locally free $\O_X$--module.
Prominent examples of free divisors are normal crossing divisors or discriminants of stable mappings $\xymat{f:X\ar[r]&Y}$ where $\dim X\ge\dim Y$.

A holomorphic function $f\in\O_{X,x}$ is called Euler homogeneous if $\chi(f)=f$ for some $\chi\in\Der_{X,x}$ which is then called an Euler vector field for $f$.
We call $f$ strongly Euler homogeneous at $x$ if it admits an Euler vector field $\chi\in\m_{X,x}\cdot\Der_{X,x}$.
The divisor $D$ is called Euler homogeneous if, for all $x\in D$, $\II(D,x)=\langle f\rangle$ for some Euler homogeneous $f\in\O_{X,x}$.
We call $D$ strongly Euler homogeneous at $x$ if $\II(D,x)=\langle f\rangle$ for some strongly Euler homogeneous $f\in\O_{X,x}$ at $x$.
By Lemma \ref{7}, a product $(D,x)=(D',x')\times(\C,0)$ with a smooth factor is Euler homogeneous and strong Euler homogeneity of $(D,x)$ and $(D',x')$ are equivalent.
By strong Euler homogeneity of $D$ we mean strong Euler homogeneity at $x$ for all $x\in D$.
This definition can be fulfilled by a single local Euler vector field only at smooth points.
Whereas Euler homogeneity is obviously an open condition, strong Euler homogeneity is not.
For example, the free divisor $D=\VV(z(x^4+xy^4+y^5))\subseteq\C^3$ \cite[Sec.~1]{CU05} is strongly Euler homogeneous at $0$, with Euler vector field $\chi=z\partial_z$, but not at the other points of the $z$--axis by Lemma \ref{7} and K.~Saito \cite{Sai71}.

The divisor $D$ is called locally quasi--homogeneous if, for all $x\in D$, $(D,x)$ is defined by a quasi--homogeneous polynomial with respect to strictly positive weights in some local coordinate system centered at $x$.
Local quasi--homogeneity obviously implies (strong) Euler homogeneity.
By K.~Saito \cite{Sai71}, the three properties are equivalent if $D$ has only isolated singularities or, in particular, in dimension $n=2$.

By F.J.~Castro--Jim\'enez et al.\ \cite{CNM96}, local quasi--homogeneity implies LCT for free divisors.
F.J.~Calder\'on--Moreno et al.\ \cite[Thm.~1.3]{CMNC02} prove that equivalence holds in dimension $n=2$ by an explicit technical construction of an Euler vector field \cite[Thm.~3.3]{CMNC02}.
In \cite[Conj.~1.4]{CMNC02}, they propose the following generalization of their result to higher dimensions which is the main motivation for this article.

\begin{cnj}[(F.J.~Calder\'on--Moreno et al.)]\label{18}
Let $D$ be a free divisor in a complex manifold $X$.
If the logarithmic comparison theorem holds for $D$ then $D$ is strongly Euler homogeneous.
\end{cnj}

\begin{thm}[(F.J.~Calder\'on--Moreno et al.)]\label{31}
Conjecture \ref{18} holds in dimension $n=2$.
\end{thm}

For any dimension $n\ge3$, F.J.~Castro--Jim\'enez and J.M.~Ucha--Enr\'iquez \cite{CU05} found a family of Euler homogeneous free divisors for which LCT does not hold.
But these divisors are not strongly Euler homogeneous since any Euler vector field at a point in $D$ with $x_1=0=x_2$ and $x_n\ne0$ has a non vanishing $\partial_n$ component.
Thus the converse of Conjecture \ref{18} is an open problem as well.

K.~Saito \cite[\S3]{Sai80} constructed the logarithmic stratification of a divisor by integration along logarithmic vector fields.
At each point of a logarithmic stratum, the logarithmic vector fields span the tangent space of this stratum.
In his language, a divisor may be called holonomic if this stratification is locally finite or, equivalently, its logarithmic characteristic variety is of minimal dimension $n$ \cite[Prop.~3.18]{Sai80}.
A divisor $D$ is called Koszul free at $x$ if an $\O_{X,x}$--basis of $\Der_{X,x}(-\log D)$ defines a regular sequence in $\gr^F\D_{X,x}$ where $F_\bullet$ denotes the filtration by order on the ring of differential operators $\D_X$ \cite[Def.~1.6]{CN02}.
Koszul free is the same as free holonomic in the above sense \cite[Cor.~1.9]{CN02}.
For example, the free divisor $D=\VV(xy(x+y)(xz+y))\subseteq\C^3$ \cite[Ex.~6.2]{CN02} fulfills LCT but is not Koszul free since any point of the $z$--axis is a logarithmic stratum.

Conjecture \ref{18} is the missing piece in the following diagram of known implications and non--implications for a free divisor \cite[Sec.~6]{CN02}.
\[
\xymat@C=10pt{
\text{locally quasi--homogeneous}\ar@{x=>}[rr]\ar@{x=>}[d]\ar@{x=>}[rrd] && \text{Euler homogeneous}\ar@{=x}[d]\\
\text{Koszul free}\ar@{x=x}[rr]\ar@{x=x}[rru] && \text{LCT}
}
\]

F.J.~Calder\'on--Moreno and L.~Narv\'aez--Macarro \cite[Cor.~4.3]{CN05} recently gave the following characterization of LCT for free divisors.

\begin{thm}[(F.J.~Calder\'on--Moreno, L.~Narv\'aez--Macarro)]\label{19}
Let $D$ be a free divisor in a complex manifold $X$.
Then LCT holds for $D$ if and only if the logarithmic Spencer complex
\[
\D_X\overset{\LL}{\otimes}_{\D_X(\log D)}\O_X(D)\simeq\D_X\otimes_{\D_X(\log D)}\Sp^\bullet_{\D_X(\log D)}(\O_X(D))
\]
is concentrated in degree $0$ and the natural morphism
\[
\xymat{
\D_X\otimes_{\D_X(\log D)}\O_X(D)\ar[r]^-{\epsilon_D}&\O_X(*D)
}
\]
is injective.
\end{thm}

The injectivity of $\epsilon_D$ in Theorem \ref{19} at $x\in D$ where $\II(D,x)=\langle f\rangle$ and $f\in\m_{X,x}$ means exactly that the annihilator $\Ann_{\D_{X,x}}(1/f)$ is generated by differential operators of order one.
For Koszul free divisors, the first condition in Theorem \ref{19} is fulfilled and the second one was characterized by T.~Torrelli \cite[Thm.~1.7]{Tor04} as follows.

\begin{thm}[(T.~Torrelli)]\label{20}
Let $0\ne f\in\m_{X,x}$ be a Koszul free germ.
Then $\Ann_{\D_{X,x}}(1/f)$ is generated by differential operators of order one if and only if $f$ is Euler homogeneous, $-1$ is the smallest integer root of the Bernstein polynomial of $f$, and the annihilator $\Ann_{\D_{X,x}}(f^s)$ is generated by differential operators of order one.
\end{thm}

In particular, LCT implies Euler homogeneity for Koszul free divisors.
We shall need only this implication but one can easily deduce the following stronger statement.

\begin{cor}\label{32}
Conjecture \ref{18} holds for Koszul free divisors.
\end{cor}

\begin{proof}
We may identify $(X,x)=(\C^n,0)$ and assume that LCT holds for the Koszul free divisor $(D,0)=\VV(f)\subseteq(\C^n,0)$ where $f\in\O_{\C^n,0}$.
By Theorems \ref{19} and \ref{20}, we know that $(D,0)$ is already Euler homogeneous.
If $(D,0)$ is not strongly Euler homogeneous then we may assume that there is an Euler vector field $\chi\in\Der_{\C^n,0}\backslash\m_{\C^n,0}\cdot\Der_{\C^n,0}$ for $f$.
By K.~Saito \cite[Lem.~3.5]{Sai80}, integration along $\chi$ yields a coordinate system $\ul x=(x_1,\dots,x_n)$ such that $f=u\cdot f'$ where $u\in\O_{\C^n,0}^*$ and $f=f(x_1,\dots,x_{n-1})\in\O_{\C^{n-1},0}$.
This means that $(D,0)=(D,0)'\times(\C,0)$ where $(D',0)=\VV(f')\subseteq(\C^{n-1},0)$ and Koszul freeness, strong Euler homogeneity, and LCT for $D$ and $D'$ are equivalent by Lemmata \ref{26} and \ref{21}.
But also the condition on the annihilator in Theorem \ref{20} is equivalent for $f$ and $f'$.
This is a contradiction by induction on the dimension $n$.
\end{proof}

In this article, we describe the formal structure of the logarithmic vector fields, that is the $\m_{X,x}$--adic completion of $\Der_{X,x}(-\log D)$.
The result in Theorem \ref{6} is obtained by performing the construction of K.~Saito \cite[\S3]{Sai71} of the Poincar\'e--Dulac decomposition \cite[Ch.~3.~\S3.2]{AA88} simultaneously to a system of generators.
In Theorem \ref{30}, we combine this result with an explicit necessary condition for LCT for a free divisor due to F.J.~Calder\'on--Moreno et al.\ \cite[Sec.~2]{CMNC02} to prove our main result:

\begin{thm}\label{30}
Conjecture \ref{18} holds in dimension $n=3$.
\end{thm}

In fact it turns out that this problem is purely formal for a non Koszul free divisor.

As a further application of Theorem \ref{6}, we describe in Propositions \ref{17} and \ref{29} the formal Lie algebra structure of the logarithmic vector fields for a free divisor in dimension $n\le3$.
In Example \ref{41}, we give a counter--example in dimension $n=4$.

\begin{acknowledgements}
The second author is grateful to David Mond for helpful hints.
We like to thank Tristan Torrelli for remarks on the introduction.
\end{acknowledgements}

\section{Vector fields}

We shall denote row vectors by a lower bar and column vectors by an upper bar.
Let $\O:=\O_n:=\C\{\ul x\}$ be the ring of convergent power series in the variables $\ul x=(x_1,\dots,x_n)$ and $\m:=\m_n:=\langle\ul x\rangle$ its maximal ideal.
There are analog definitions and statements as in this section for the ring $\wh\O=\C[\![\ul x]\!]$ with maximal ideal $\wh\m$.
The $\C$--linear derivations of $\O$ form the module $\Der:=\Der_n:=\Der_\C\O_n$ of vector fields.
It is a free $\O$--module of rank $n$ with basis the partial derivatives $\ul\partial=(\partial_1,\dots,\partial_n)$.

The module $\Der$ acts naturally on $\O$ and on itself by the Lie bracket $\delta(\eta):=[\delta,\eta]$ where $\delta,\eta\in\Der$.
Weights $\ul w=(w_1,\dots,w_n)\in\C^n$ define a vector field $\sigma=\sum_iw_ix_i\partial_i\in\Der$.
A power series $p\in\O$ or a vector field $\delta\in\Der$ is $\ul w$-- or $\sigma$--homogeneous of degree $\lambda\in\C$ if $\sigma(p)=\lambda\cdot p$ or $\sigma(\delta)=[\sigma,\delta]=\lambda\cdot\delta$.
When referring to the standard weights $\ul w=(1,\dots,1)$, we omit $\ul w$.
In this case, the $x_i$ are homogeneous of degree $1$, the $\partial_i$ of degree $-1$.

\begin{ntn} 
Any vector field $\delta\in\Der$ can be uniquely written as $\delta=\sum_{i=-1}^\infty\delta_i$ where $\delta_i$ is homogeneous of degree $i$ and $\delta_0=\ul x A\ol\partial$ for a unique matrix $A\in\C^{n\times n}$.
For $\delta\in\m\cdot\Der$, we call $\delta_0$ the linear part of $\delta$.
\end{ntn}

\begin{lem}\label{8}
$[\ul x A\ol\partial,\ul x B\ol\partial]=\ul x[A,B]\ol\partial$.
\end{lem}

\begin{proof}
This follows immediately from $\ol\partial(\ul x)=(\partial_i(x_j))_{i,j}=(\delta_{i,j})_{i,j}$.
\end{proof}

\begin{dfn}
Let $\delta= \sum_{i=0}^\infty\delta_i\in\m\cdot\Der$ and let $A\in\C^{n\times n}$ be such that $\delta_0=\ul x A\ol\partial$.
Then $\delta$ is called semisimple (resp.\ diagonal) if $\delta=\delta_0$ and $A$ is semisimple (resp.\ diagonal).
It is called nilpotent if $A$ is nilpotent.
\end{dfn}

\begin{lem}\label{10}
A nilpotent $\delta\in\m\cdot\Der$ is nilpotent on $\O/\m^k$ and 
$\Der/\m^k\cdot\Der$ for all $k\ge 0$.
\end{lem}

\begin{proof}
Let $\delta\in\m\cdot\Der$ where $\delta_0=\ul x A\ol\partial$ and $A\in\C^{n\times n}$ is nilpotent.
After a $\C$--linear coordinate change, we may assume that $A$ has Jordan normal form.
Order the monomials first by minimal degree in $\ul x$ and then lexicographically by $\partial_1>\cdots>\partial_n>x_n>\cdots>x_1$.
Then 
\begin{align*}
x_i\partial_{i+1}(x_1^{\alpha_1}\cdots x_n^{\alpha_n})&=\alpha_{i+1}x_1^{\alpha_1}\cdots x_i^{\alpha_i+1}x_{i+1}^{\alpha_{i+1}-1}\cdots x_n^{\alpha_n}<x_1^{\alpha_1}\cdots x_n^{\alpha_n},\\
x_i\partial_{i+1}(\ul x^{\ul \alpha}\partial_j)&=x_i\partial_{i+1}(\ul x^{\ul \alpha})\partial_j-\delta_{i,j}\ul x^{\ul \alpha}\partial_{i+1},
\end{align*}
and $x_i\partial_{i+1}(\ul x^{\ul \alpha})\partial_j,\ul x^{\ul \alpha}\partial_{j+1}<\ul x^{\ul \alpha}\partial_j$ which implies the claim.
\end{proof}

Recall that $\delta_0$ is the linear part of $\delta\in\m\cdot\Der.$

\begin{ntn}
Any vector field $\delta\in\m\cdot\Der$ can be uniquely written as $\delta=\delta_S+\delta_N$ where $\delta_S$ is semisimple, $\delta_N$ is nilpotent, and $[\delta_{S,0},\delta_{N,0}]=0$. 
Note that $\delta_S=\delta_{S,0}$ and $\delta_N=\delta_{N,0}+\sum_{i=1}^\infty\delta_i$.
\end{ntn}

Note that semisimplicity (resp.\ nilpotency) of $\delta\in\m\cdot\Der$ means that $\delta = \delta_S$ (resp.\ $\delta = \delta _N$).
By Lemma \ref{8}, for weights $\ul w\in\C^n$ and a $\ul w$--homogeneous $\delta\in\m\cdot\Der$ of
degree $0$, there is a $\ul w$--homogeneous $\C$--linear coordinate 
change such that $\delta_0$ is in Jordan normal form.
In particular, one can always assume that $\delta_S$ is diagonal in this case.

\begin{lem}\label{9}
Let $\ul w=(w_1,\dots,w_n)\in\Q^n$ be rational weights and let $\delta\in\m\cdot\Der$ be a $\ul w$--homogeneous vector field of degree $\lambda\in\Q^*$.
Then $\delta$ is nilpotent.
\end{lem}

\begin{proof}
Let $\sigma:=\sum_iw_ix_i\partial_i=:\ul x D\ol\partial\in\Der$ where $D\in\Q^{n\times n}$ is diagonal and $A=(a_{i,j})_{i,j}\in\C^{n\times n}$ such that $\delta_0=\ul x A\ol\partial$.
By Lemma \ref{8}, $\sigma(\delta)=\lambda\cdot\delta$ implies that
\[
\bigl((w_i-w_j)\cdot a_{i,j}\bigr)_{i,j}
=[D,A]=
\lambda\cdot A=(\lambda\cdot a_{i,j})_{i,j}.
\]
We may assume that $\lambda>0$ and $w_1\ge\cdots\ge w_n$.
But then $a_{i,j}=0$ for $i\le j$ and hence $\delta$ is nilpotent.
\end{proof}

Nilpotency of vector fields is clearly invariant under arbitrary coordinate changes.
We shall see that diagonal vector fields are invariant under coordinate changes which are homogeneous for the corresponding weights.

\begin{lem}\label{5}
Let $\sigma=\sum_iw_ix_i\partial_i\in\Der$ and $\ul w:=(w_1,\dots,w_n)\in\C^n$.
Then $\sigma$ is invariant under $\ul w$--homogeneous coordinate changes.
\end{lem}

\begin{proof}
Let $y_i=x_i+h_i$ with $\sigma(h_i)=w_ih_i$.
Then $\partial_{x_i}=\partial_{y_i}+\sum_j\frac{\partial h_j}{\partial x_i}\partial_{y_j}$ and hence
\begin{align*}
\sum_iw_iy_i\partial_{y_i}
&=\sum_iw_i(x_i+h_i)\Bigl(\partial_{x_i}-\sum_j\frac{\partial h_j}{\partial x_i}\partial_{y_j}\Bigr)\\
&=\sum_iw_ix_i\partial_{x_i}+w_ih_i\partial_{x_i}\\
&-\sum_{i,j}w_ix_i\frac{\partial h_j}{\partial x_i}\partial_{y_j}+w_ih_i\frac{\partial h_j}{\partial x_i}\partial_{y_j}\\
&=\sum_iw_ix_i\partial_{x_i}+w_ih_i\partial_{x_i}\\
&-\sum_jw_jh_j\partial_{y_j}-\sum_iw_ih_i(\partial_{x_i}-\partial_{y_i})=\sum_iw_ix_i\partial_{x_i}.
\end{align*}
\end{proof}

\section{Logarithmic vector fields}

Let $0\ne f\in\O$ be a convergent power series.
There are analog definitions and statements as in this section for a formal power series $0\ne f\in\wh\O$.

\begin{dfn}\label{25}
The $\O$--module of logarithmic vector fields is defined by
\[
\Der_f:=\Der_f\O:=\{\delta\in\Der\mid\delta(\O\cdot f)\subseteq\O\cdot f\}.
\]
If $\Der_f\not\subseteq\m\cdot\Der$ then we call $f$ a product (with a smooth factor).
If $\chi(f)=f$ for some $\chi\in\Der_f$ then $f$ is called Euler homogeneous and $\chi$ is called an Euler vector field (for $f$).
We call $f$ strongly Euler homogeneous (at the origin) and $\chi$ a strong Euler vector field (for $f$) if $\chi\in\m\cdot\Der$.
\end{dfn}

$\Der_f$ is invariant under multiplication of $f$ by units and the Lie bracket on $\Der$ induces a Lie bracket on $\Der_f$.
By the Leibniz rule,
\[
\Der_f=\{\delta\in\Der\mid \delta(f)\in\O\cdot f\}
\]
can be identified with the projection of the first syzygy module of 
$\frac{\partial f}{\partial x_1},\dots,\frac{\partial f}{\partial
 x_n},f$ to the first $n$ components.
In particular,
\[
\Der_f\wh\O=\wh{\Der_f\O}
\]
is the $\m$--adic completion of $\Der_f\O$ and $f$ being a product is invariant under completion.
Euler homogeneity of $f$ is equivalent to $\delta(f)\notin\m\cdot f$ for some $\delta\in\Der_f$.
In particular, Euler homogeneity is invariant under completion and strong Euler homogeneity at the origin as well.
Moreover, strong Euler homogeneity of $f$ is invariant under multiplication of $f$ by units.
Indeed, if $\chi\in\Der_f$ is a strong Euler vector field for $f$ and $u\in\O^*$ then $u(u+\chi(u))^{-1}\chi$ is a strong Euler vector field for $u\cdot f$.

\begin{lem}\label{7}
If $f\in\O_n$ is a product then there is a coordinate change such that $f=u\cdot f'$ for some unit $u\in\O_n^*$ and some $f'\in\O_{n-1}$ .
In this case, 
\[
\Der_f\O_n=\O_n\cdot\Der_{f'}\O_{n-1}\oplus\O_n\cdot\partial_n,
\]
$\partial_n$ is an Euler vector field for $\exp(x_n)\cdot f'$, and strong Euler homogeneity of $f$ and $f'$ are equivalent.
\end{lem}

\begin{proof}
A more general version of the first statement is given by K.~Saito \cite[Lem.~3.5]{Sai80}.
If $\chi=\chi'+a_n\partial_n\in\m\cdot\Der$ where $\chi'\in\m_n\cdot\Der_{n-1}$ and $a_n\in\m$ is an Euler vector field for $f'$ then $\chi'_{\vert x_n=0}\in\m_{n-1}\cdot\Der_{n-1}$ is also an Euler vector field for $f'$ which implies the last statement.
\end{proof}

\begin{lem}
Strong Euler vector fields are non nilpotent.
\end{lem}

\begin{proof}
Choose $k$ such that $f\in\m^k\backslash\m^{k+1}$.
Then $[f]\in\m^k/\m^{k+1}$ is an eigenvector with eigenvalue $1$ of any strong Euler vector field.
Therefore such a vector field can not be nilpotent by Lemma \ref{10}.
\end{proof}

\begin{lem}\label{39}
$\Der_{g\cdot h}=\Der_g\cap\Der_h$ for all $f,g\in\O$.
\end{lem}

\begin{proof}
Let $f=f_1^{l_1}\cdots f_m^{l_m}$ be a decomposition of $f\in\O$ into irreducible factors.
Then 
\[
\delta(f)=\sum_{i=1}^ml_i\cdot\delta(f_i)\cdot f/f_i,
\]
for all $\delta\in\Der$, and hence $\Der_f=\bigcap_{i=1}^m\Der_{f_i}$.
\end{proof}

\section{Freeness and Saito's Criterion}

\begin{dfn}\label{15}
We call a reduced $f\in\O$ free if $\det\bigl(\delta_i(x_j)\bigr)_{i,j}\in\O^*\cdot f$ for some elements $\ul\delta=(\delta_1,\dots,\delta_n)\in\Der_f$.
Note that $\ol\delta=A\ol\partial$ for $A:=\bigl(\delta_i(x_j)\bigr)_{i,j}\in\O^{n\times n}$ and $\ol\delta:=(\delta_1,\dots,\delta_n)\in\Der_f$.
We define freeness of $f\in\wh\O$ analogously.
\end{dfn}

Freeness of $f$ is invariant under coordinate changes, multiplication of $f$ by a unit, and completion.
By Saito's Criterion \cite[Lem.~1.8.ii]{Sai80}, a convergent $f\in\O$ is free if and only if $\Der_f$ is a free $\O$--module and $\ul\delta$ in Definition \ref{15} is an $\O$--basis of $\Der_f$.
One of these implications also holds for a formal $f\in\wh\O$.

\begin{prp}[(Formal Saito's Criterion)]
If $f\in\wh\O$ is free then $\Der_f$ is a free $\wh\O$--module of rank $n$ and $\ul\delta$ in Definition \ref{15} is an $\wh\O$--basis of $\Der_f$.
\end{prp}

\begin{proof}
The statement is obvious for $f\in\wh\O^*$ and we may assume that $f\in\wh\m$.
We first prove that $f\mid\det A$ for any $\ol\delta\in\Der_f^n$ and $A$ as in Definition \ref{15}. 
Let $\ol a:=\ol\delta(f)/f\in\wh\O^n$, $B\in\wh\O^{n\times n}$ the adjoint matrix of $A$, and $f=f_1\cdots f_k$ a decomposition of $f$ into different irreducible factors $f_i\in\wh\m$.
Then $\det A\cdot\ol\partial(f)=B\ol a\cdot f$ and hence
\[
f\mid\sum_{i=1}^k\det A\cdot f_1\cdots\partial_j(f_i)\cdots f_k
\]
for all $j=1,\dots,n$.
Then $f_i\mid\det A\cdot f_1\cdots\partial_j(f_i)\cdots f_k$ for all $i=1,\dots,k$ and $j=1,\dots,n$.
For some $j$, $f_i\mid\partial_j(f_i)$ is impossible and hence $f_i\mid\det A$ for all $i=1,\dots,k$ and finally $f\mid\det A$.

Now assume that $\det A\in\wh\O^*\cdot f$.
Then $\ul\delta$ is $\wh\O$--linearly independent.
Let $\delta\in\Der_f$ and $\ol b\in\wh\O^n$ such that $\delta=\ul b\ol\partial$.
Then $f\mid\ul bB$ by the preceding arguments and hence
\[
\delta=\ul b\ol\partial=f^{-1}\ul bBA\ol\partial=f^{-1}\ul bB\ol\delta\in\langle\ul\delta\rangle.
\]
Thus $\ul\delta$ generates $\Der_f$ and is an $\wh\O$--basis.
\end{proof}

\section{Formal Structure Theorem}

This section concerns only formal power series.
The results of K.~Saito in \cite[\S\S2-3]{Sai71} are compatible with multiweights $W=(\ul w^1,\dots,\ul w^s)$ where $\ul w^i=(w^i_1,\dots,w^i_n)\in\C^n$ in the following sense.
 
\begin{lem}[(K.~Saito {\cite[Lem.~2.3.iii]{Sai71}})]\label{2}
Any $p\in\wh\O$ can be uniquely written as $p=\sum_{\ul\lambda\in\C^s}p_{\ul\lambda}$ where $p_{\ul\lambda}\in\wh\O$ is $W$--multihomogeneous of degree $\ul\lambda$.
\end{lem}

\begin{lem}[(K.~Saito {\cite[Kor.~2.5]{Sai71}})]\label{3}
Let $\delta=\delta_0\in\Der$ be linear and $W$--multihomogeneous of degree $\ul0$. 
Assume that its semisimple part $\delta_S=\sum_iw_ix_i\partial_i$ is diagonal and set $\ul w:=(w_1,\dots,w_n)\in\C^n$.
Then, for any $W$--multihomogeneous $p\in\wh\O$ of degree $\ul\lambda\in\C^s$ and any $\lambda\in\C$, there is a $W$--multihomogeneous $q\in\wh\O$ of degree $\ul\lambda$ such that $\delta(q)-\lambda\cdot q+p$ is $\ul w$--homogeneous of degree $\lambda$.
\end{lem}

\begin{proof}
In the proof of \cite[Lem.~2.4]{Sai71}, we only need to replace the space $W_{m,\lambda}$ of homogeneous polynomials of degree $m$ which are $\ul{w}$--homogeneous of degree $\lambda$ by the subspace $W_{m,\ul{\lambda },\lambda }\subseteq W_{m,\lambda }$ of $W$--multihomogeneous elements of degree $\ul{\lambda }$. 
Since this space is stable by $\delta$, the same linear algebra argument applies and the claim follows exactly in the same way as in \cite[Kor.~2.5]{Sai71}.
\end{proof}

\begin{thm}[(K.~Saito {\cite[Satz 3.1]{Sai71}})]\label{4}
Let $\delta\in\wh\m\cdot\Der$ be $W$--multihomogeneous of degree $\ul0$, $\delta_S=\sum_iw_ix_i\partial_i$, and $\ul w:=(w_1,\dots,w_n)\in\C^n$.
Then $\delta$ is $\ul w$--homogeneous of degree $0$ after a $W$--homogeneous coordinate change.
In particular, $[\delta_S,\delta_N]=0$ in this case.
\end{thm}

\begin{proof}
The same proof as in \cite[Satz 3.1]{Sai71} works with a sequence of coordinate changes $x_i^{(m)}=x_i^{(m-1)}+h_i$ tangent to the identity.
We only need to add the condition that each $h_i$ is $W$--multihomogeneous of the same degree as $x_i$.
In fact, this follows from the recursion formulas used by K.~Saito, precisely because of Lemma \ref{3}, with the multihomogeneity of the coefficients in
\[
\delta=\sum g_i^{(m)}\bigl(\ul{x} ^{(m)}\bigr)\cdot \frac{\partial }{\partial x_i^{(m)}}
\]
proved simultaneously.
\end{proof}

The following result is a formal structure theorem for $\Der_f$.
For reduced convergent $f$, $\Der_f$ depends only on the zero set of $f$ or the divisor defined by $f$.
Considering this divisor means to consider $f$ up to contact equivalence which allows coordinate changes and multiplication of $f$ by units.
But invariance of $\Der_f$ under contact equivalence also holds for non reduced and formal $f$.
We define the formal divisor $\wh\VV(f)$ associated with $f$ as the formal contact equivalence class of $f$.
The invariant $s$ defined below can be considered as the maximal multihomogeneity of an equation of this formal divisor. 

Let $0\ne f\in\wh\O$ be a formal power series. 
We assume that $f$, considered as $\wh\VV(f)$, is not a product which means, by definition, that $\Der_f\subseteq\m\cdot\Der$.

\begin{thm}[(Formal Structure Theorem)]\label{6}
Let $s$ be the maximal dimension of the vector space of diagonal $\sigma\in\Der_f$ with $\sigma(f)\in\C\cdot f$, for $f$ varying in a formal contact equivalence class.
This means that $s$ is maximal for all coordinate systems and changes of $f$ by a factor in $\wh\O^*$.
Then there are $\sigma_1,\dots,\sigma_s,\nu_1,\dots,\nu_r\in\Der_f$, a coordinate change, a change of $f$ by a factor in $\wh\O^*$, and a set of irreducible factors $f_1,\dots,f_m$ of $f$ such that
\begin{asparaenum}[1.]
\item $\sigma_1,\dots,\sigma_s,\nu_1,\dots,\nu_r$ is a minimal system of generators of $\Der_f$,
\item\label{38} if $\delta\in\Der_f$ with $[\sigma_i,\delta]=0$ for all $i$ then $\delta_S\in\langle\sigma_1,\dots,\sigma_s\rangle_\C$, 
\item $\sigma_i$ is diagonal with eigenvalues in $\Q$,
\item $\nu_i$ is nilpotent,
\item $[\sigma_i,\nu_j]\in\Q\cdot\nu_j$, and
\item $\sigma_i(f_j)\in\Q\cdot f_j$.
\end{asparaenum}
\end{thm}

\begin{proof}
Let $\sigma_1,\dots,\sigma_s\in\Der_f$ where $\sigma_i=\sum_jw^i_jx_j\partial_j$ and $w^i_j\in\C$ such that $\sigma_i(f)\in\C\cdot f$.
By Lemma \cite[Lem.~1.4]{Sai71}, we may assume that $\ul w^i:=(w^i_1,\dots,w^i_n)\in\Q^n$ and we denote $W:=(\ul w^1,\dots,\ul w^s)$.
Then $f$ is $W$--multihomogeneous of some degree $\ul\lambda=(\lambda_1,\dots,\lambda_s)\in\Q^s$.
Let $\delta\in\Der_f\backslash\langle\sigma_1,\dots,\sigma_s\rangle_\C$ and $a\in\wh\O$ such that $\delta(f)=a\cdot f$.
By Lemma \ref{2}, we may assume that $\delta$ and $a$ are $W$--multihomogeneous of some degree $\ul\mu\in\Q^s$.
By Lemma \ref{9}, $\delta$ is nilpotent if $\ul\mu\ne\ul0$ and we may hence assume that $\ul\mu=\ul0$.
By Theorem \ref{4}, there is a $W$--multihomogeneous coordinate change such that $\delta$ is $\ul w$--homogeneous of degree $0$ where $\delta_S=\sum_iw_ix_i\partial_i$ and $\ul w:=(w_1,\dots,w_n)\in\C^n$.
By Lemma \ref{5}, the $\sigma_i$ are invariant under this coordinate change.

We shall multiply $f$ by a $W$--multihomogeneous $u\in\wh\O^*$ of degree $\ul 0$ to make $a$ $\ul w$--homogeneous of degree $0$.
The transformation of $a$ under this operation is given by
\[
\delta(uf)=\bigl(\delta(u)\cdot u^{-1}+a\bigr)\cdot uf.
\]
Let $d\ge1$ be the minimal degree in which $a$ is not $\ul w$--homogeneous of degree $0$.
If $u=1+u_d$ where $u_d$ is homogeneous of degree $d$ then $\delta_0(u_d)+a_d$ is the degree $d$ part of the transformed $a$.
By Lemma \ref{3}, there is a $W$--multihomogeneous $u_d$ of degree $\ul 0$ such that $\delta_0(u_d)+a_d$ is $\ul w$--homogeneous of degree $0$.
Then the desired $u$ exists by induction on $d$.
We may hence assume that $a$ is $\ul w$--homogeneous of degree $0$.

Let $f=\sum_\lambda f_\lambda$ be the expansion of $f$ in $\ul w$--homogeneous parts as in Lemma \ref{2}.
Then $\delta(f_\lambda)=a\cdot f_\lambda$ and hence $\delta_N(f_\lambda)=(a-\lambda)\cdot f_\lambda$ for all $\lambda$.
By Lemma \ref{10}, $f_\lambda=0$ for $\lambda\ne a_0$ and hence $f=f_{a_0}$.
But then $\delta_S(f)=a_0\cdot f$ and hence, by the minimality assumption, $\delta_S\in\langle\sigma_1,\dots,\sigma_s\rangle_\C$.
So we can assume that $\delta=\delta_N$ is nilpotent.

It remains to construct the $W$--homogeneous irreducible factors of $f$.
Let $W':=(\ul w^1,\dots,\ul w^{s'})$ where $s'\le s$ is maximal such that there is a set of irreducible $W'$--homogeneous factors $f_1,\dots,f_m$ of $f$.
We may assume that $m\ge2$ and $s'<s$.
Let $t:=s'+1$, $\sigma_t=\sum_iw_ix_i\partial_i$, $\ul w:=(w_1,\dots,w_n)\in\Q^n$, and $f=f_1^{l_1}\cdots f_m^{l_m}$.
By Lemma \ref{39}, $\sigma_t(f_i)=a_{t,i}\cdot f_i$ for some $a_{t,i}\in\wh\O$.
By the above argument, there are, for $i=1,\dots,m-1$, $W'$--multihomogeneous $u_{t,i}\in\wh\O^*$ and $\lambda_{t,i}\in\Q$ such that 
\begin{equation}\label{40}
\sigma_t(u_{t,i}f_i)=\lambda_{t,i}\cdot u_{t,i}f_i.
\end{equation}
We choose $u_{t,m}\in\wh\O^*$ such that $\prod_{i=1}^mu_{t,i}^{l_i}=1$ and set $\lambda_{t,m}:=l_m^{-1}\cdot\bigl(\lambda_t-\sum_{i=1}^{m-1}l_i\lambda_{t,i}\bigr)\in\Q$.
Then (\ref{40}) holds for $i=1,\dots,m$ and hence $u_{t,1}f_1,\dots,u_{t,m}f_m$ form a set of $W'$--multihomogeneous, $\ul w$--homogeneous, irreducible factors of $f$.
This contradicts to the maximality of $s'$ and finishes the proof.
\end{proof}

\begin{cor}\label{22}
If $f$ is Euler homogeneous then at least one $\sigma_i$ in Theorem \ref{6} can be chosen to be an Euler vector field.
\end{cor}

\begin{proof}
This follows immediately from Lemma \ref{10}.
\end{proof}

\begin{cor}\label{12}
If $\Der_f$ is a free $\wh\O$--module then there is a basis of $\Der_f$ as in Theorem \ref{6}.
\end{cor}

\begin{proof}
This follows immediately from Nakayama's Lemma.
\end{proof}

\begin{cor}\label{16}
In Theorem \ref{6}, if $f$ is free and $\lambda^i_j\in\Q$ such that $[\sigma_i,\nu_j]=\lambda^i_j\cdot\nu_j$ then $f$ is $\sigma_i$--homogeneous of degree $\sum_{j=1}^nw^i_j+\sum_{j=1}^{n-s}\lambda^i_j$.
\end{cor}

\begin{proof}
This follows immediately from Definition \ref{15}.
\end{proof}

\section{Formal Lie algebra structure}

This section concerns only formal power series.
Let $0\ne f\in\wh\O$ be not a product and let $\sigma_1,\dots,\sigma_s,\nu_1,\dots,\nu_r\in\Der_f$ be as
in Theorem \ref{6}.
Let $\DD_d$ be the Lie algebra $\Der_f/\m^d\cdot\Der_f$ over $\C$ where $d\ge1$.

In this section, we shall freely denote by the same letter a
vector field $\delta\in\Der_f$ and its class modulo $\m^d\cdot\Der_f$, $\delta
\in\DD_d$.
Then $\SS_d:=\bigoplus_{i=1}^s\C\cdot\sigma_i\subseteq\DD_d$ is an abelian Lie subalgebra.
The centralizer $\Cen(\SS_d)$ of $\SS_d$ in $\DD_d$ is the Lie subalgebra of $\ul\sigma$--multihomogeneous logarithmic vector fields of degree $\ul0$ where $\ul\sigma:=(\sigma_1,\dots,\sigma_s)$.
By Theorem \ref{6}, $\Cen(\SS_d)$ is an almost algebraic Lie algebra \cite[III.11]{Jac62}.
The derived series of $\DD_d$ is defined by
\begin{alignat*}{2}
\DD_d^{(0)}&:=\DD_d, &\quad \DD_d^{(i+1)}&=[\DD_d^{(i)},\DD_d^{(i)}]
\end{alignat*}
and $\DD_d$ is called solvable if $\DD_d^{(i)}=0$ for $i\gg0$ \cite[I.7]{Jac62}.
By Lemma \ref{10} and Engel's Theorem \cite[II.3]{Jac62},
$\m\cdot\DD_d$ is a nilpotent ideal and hence $\DD_d$ is solvable if
and only if $\DD_d/\m\cdot\DD_d=\DD_1$ is solvable.
An element $\delta\in\DD_1$ is reduced to its linear part
$\delta=\ul xA\ol \delta$ where $A\in\C^{n\times n}$.
In the rest of this section, we shall hence assimilate $\delta$ to the matrix $A$.

\begin{prp}\label{17}
The Lie algebras $\DD_d$ are solvable if $r\le1$ or $s=0$ in Theorem \ref{6}.
In particular, this holds if $\Der_f$ is a free $\wh\O$--module of rank $2$.
\end{prp}

\begin{proof}
If $r=0$ then $\DD_1=\SS_1$ is abelian.
If $r=1$ then $\DD_1=\SS_1\oplus\C\cdot\nu_1$ and $\C\cdot\nu_1\subseteq\DD_1$ is a nilpotent ideal.
If $s=0$ then $\DD_1$ is nilpotent by Theorem \ref{6}.\ref{38} and Engel's Theorem \cite[II.3]{Jac62}.
The second claim follows from Corollary \ref{12}.
\end{proof}

In the following, we shall prove the solvability of $\DD_1$ for a free $f$ in dimension $n=3$.
By Proposition \ref{17}, it suffices to consider the case $s=1$ in Theorem \ref{6}.
In a convenient system of coordinates $(x,y,z)$, $\Der_f$ is generated by a diagonal vector field $\sigma = ax\partial_x+by\partial_y+cz\partial_z$ where $a,b,c\in\Q$ and two $\sigma$--homogeneous nilpotent vector fields $\nu_1,\nu_2$ of degrees $\lambda_1,\lambda_2$.
The set of eigenvalues of $[\sigma,\cdot]$ is the set of differences of $a,b,c$ and includes $\lambda_1$ and $\lambda_2$.
There is $\sigma$--homogeneous relation
\[
\mu:=[\nu_1,\nu_2]=\lambda\sigma+p_1\nu_1+p_2\nu_2\quad\text{ where }\quad\lambda,p_1,p_2 \in\C
\]
of degree $\lambda_1+\lambda_2$.
As the trace of a commutator, $\tr\mu=0$ and, by additivity, $\tr(\lambda\sigma)=0$. 

We first show that $\lambda \neq 0$ if the Lie algebra $\DD_1$ is not solvable.
Indeed if $\lambda = 0$ then $\mu=p_1\nu_1+p_2\nu_2$ and $\DD_1^{(1)}\subseteq\langle\nu_1,\nu_2\rangle$. 
Therefore $\DD_1^{(2)}\subseteq\C\cdot\mu$ and hence $\DD_1^{(3)}=0$ which proves that $\DD_1$ is solvable.
(In fact we might prove with some more calculations that already $\mu=0$.)

We may assume now that $\lambda\neq 0$.
Then the $\sigma$--degree of $\lambda\sigma$ and hence of $p_1\nu_1$, $p_2\nu_2$, and $\mu$ equals $0$.
In particular, $\lambda_i\ne0$ implies $p_i=0$ for $i=1,2$ and $\lambda_1+\lambda_2=0$ being the $\sigma$--degree of $\mu$.
Finally the situation of a non solvable $\DD_1$ reduces to the following two cases: 
\begin{asparaenum}[{\bf Case} I:]
\item\label{34} $\mu=\lambda\sigma+p_1\nu_1+p_2\nu_2$ where $\lambda \neq 0$ and $\lambda_1=\lambda_2=0$.
\item\label{35} $\mu=\lambda\sigma$ where $\lambda \neq 0$ and $\lambda_1=-\lambda_2\neq 0$.
\end{asparaenum}
In Case \ref{34}, we would have $\DD_1^{(1)}=\C\cdot\mu$ and $\DD_1$ would be solvable as follows from $\DD_1^{(2)}=0$.
However, we shall prove easily that Case \ref{34} can not occur.
We shall also prove that Case \ref{35} is impossible by a more complicate argument.

\begin{asparaenum}[{\bf Case} I:]

\item $\mu=\lambda\sigma+p_1\nu_1+p_2\nu_2$ where $\lambda\neq 0$ and $\lambda_1=\lambda_2=0$.
We may assume that $a\ne0$ and let $E$ be the $a$--eigenspace of $\sigma$.
The equality $\lambda_i=0$ means that the $\nu_i$ commute with $\sigma$ and that $E$ is invariant under the $\nu_i$ and hence under $\mu$.
By restricting to $E$, we obtain the contradiction $0=\tr(\sigma_{\vert E})=a\cdot\dim E$.

\item $\mu=\lambda\sigma$ where $\lambda \neq 0$ and $\lambda_1=-\lambda_2\neq 0$.

\begin{asparaenum}[{\em Subcase} a:] 

\item $\sigma = a(x\partial_x+y\partial_y+ z\partial_z)$.
This is impossible since all $\sigma$--homogeneous vector fields are of degree $a-a=0$.

\item $\sigma = ax\partial_x+b( y\partial_y+ z\partial_z)$ where $a\ne b$.
We may assume that $\lambda_1=a-b$ and write 
\[
\nu_1 = \begin{pmatrix}
0 & r & s \\
0 & 0 & 0 \\
0 & 0 & 0
\end{pmatrix}
\quad\text{ and }\quad 
\nu_2 = \begin{pmatrix}
0 & 0 & 0 \\
u & 0 & 0 \\
v & 0 & 0
\end{pmatrix}.
\]
Calculating the commutator gives
\[
\begin{pmatrix}
ru+sv & 0 & 0 \\
0 & -ur & -us \\
0 & -vr & -vs 
\end{pmatrix} 
= [\nu_1 , \nu_2] = \lambda\sigma =
\lambda\cdot
\begin{pmatrix}
a & 0 & 0 \\
0 & b & 0 \\
0 & 0 & b 
\end{pmatrix} 
\]
which is equivalent to
\[
vr=us=0\quad\text{ and }\quad(ru+sv,ur,vs)=(\lambda a,\lambda b,\lambda b).
\]
The equations $vr=0=us$ imply at least one of the equalities $v=0$, $r=0$, $u=0$, $s=0$. 
Each of these taken into the other relations $\lambda b= ur=vs=\frac{1}{2}\lambda a$ gives $b=0$ and then $a=0$, a contradiction.

\item
$\sigma = ax\partial_x+by\partial_y+ cz\partial_z$ where $a\ne b\ne c\ne a$.
In this case, the relations between $\sigma,\nu_1,\nu_2$ alone do not contradict to a non solvable $\DD_1$. 
We shall exclude this case by using the equation $f$. 

Up to permutation, there are two cases:
$\lambda_1=a-c$ and $\lambda_2=c-a$ where $\pm(c-b)\neq\lambda_1\neq\pm(a-b)$ or $\lambda_1=a-b$ and $\lambda_2=c-b$.

In the first case, we may assume by changing $\nu_1$, $\nu_2$, and $\sigma$ (or
equivalently $\lambda$) by a constant factor that 
\[
\nu_1 = \begin{pmatrix}
0 & 0 & 1 \\
0 & 0 & 0 \\
0 & 0 & 0
\end{pmatrix},
\quad
\nu_2 = \begin{pmatrix}
0 & 0 & 0 \\
0 & 0 & 0 \\
1 & 0 & 0
\end{pmatrix},
\quad
\sigma = \begin{pmatrix}
1 & 0 & 0 \\ 
0 & 0 & 0 \\
0 & 0 & -1
\end{pmatrix}
= \frac{1}{2} [\nu _1,\nu _2].
\]
So the linear parts of the generators of $\Der_f$ are in the canonical form
\begin{equation}\label{36}
\sigma = x\partial_x-z\partial_z, \quad \nu_1=x\partial_z, \quad \nu_2=z\partial_x.
\end{equation}
In the second case, we may assume that $(a,b,c)=(1,0,-1)$ after changing $\sigma$ by a
constant factor since $\tr\sigma=0$. 
Then we may write
\[
\nu_1 = \begin{pmatrix}
0 & 1 & 0 \\
0 & 0 & r \\
0 & 0 & 0
\end{pmatrix}
\quad\text{ and }\quad 
\nu_2 = \begin{pmatrix}
0 & 0 & 0 \\
s & 0 & 0 \\
0 & t & 0
\end{pmatrix}.
\]
Calculating the commutator gives
\[
\begin{pmatrix}
s & 0 & 0 \\
0 & rt-s & 0 \\
0 & 0 & -rt 
\end{pmatrix} 
= [\nu_1 , \nu_2] = \lambda \sigma =
\lambda\cdot
\begin{pmatrix}
1 & 0 & 0 \\
0 & 0 & 0 \\
0 & 0 & -1
\end{pmatrix}
\]
and hence $s=\lambda =rt$.
By the change of coordinates $(x,y,rz)$, we reduce to the case $r=1$ and, by dividing $\nu_2$ by $\lambda$, to the case $s=t=1$.
So the linear parts of the generators of $\Der_f$ are in the canonical form
\begin{equation}\label{37}
\sigma = x\partial_x-z\partial_z, \quad \nu_1=x\partial_y+y\partial_z, \quad \nu_2=y\partial_x+z\partial_y.
\end{equation}

By Corollary \ref{16}, we may assume that $f$ is $\sigma$--homogeneous of degree $0$ with the same $\sigma$ for (\ref{36}) and (\ref{37}).
We can hence write $f=\sum_{i\ge k}f_i$ where $f_i$ is homogeneous of degree $i$ and $0\ne f_k=\sum_ic_ix^iy^{k-2i}z^i$. 
By Lemma \ref{10}, $\nu_1(f_k)=\nu_2(f_k)=0$.
For (\ref{36}), we obtain $\frac{\partial f_k}{\partial x}=\frac{\partial f_k}{\partial z}=0$ and hence $f_k=f_k(y)$ depends only on $y$.
For (\ref{37}), the coefficients of the equation $\nu_1(f_k)=0$ are $(k-2i+2)c_{i-1}+ic_i=0$ and hence $c_0\ne0$.
Thus both (\ref{36}) and (\ref{37}) contradict to $f\in\langle x,z\rangle$ by Definition \ref{15}.

\end{asparaenum}

\end{asparaenum}

\begin{prp}\label{29}
The Lie algebras $\DD_d$ are solvable for free $f$ in dimension $n\le3$.
\end{prp}

\begin{proof}
This follows from Proposition \ref{17} and the preceding arguments.
\end{proof}

\begin{exa}\label{41}
Consider the representation of the non--solvable Lie algebra $\C\times\mathfrak{sl}_2$ defined by
\[
X=
\begin{pmatrix}
1 & 0 & 0 & 0 \\
0 & 1 & 0 & 0 \\
0 & 0 & 1 & 0 \\
0 & 0 & 0 & 1
\end{pmatrix},
\ 
H=
\begin{pmatrix}
3 & 0 & 0 & 0 \\
0 & 1 & 0 & 0 \\
0 & 0 & -1 & 0 \\
0 & 0 & 0 & -3
\end{pmatrix},
\ 
S_+=
\begin{pmatrix}
0 & 0 & 0 & 0 \\
1 & 0 & 0 & 0 \\
0 & 1 & 0 & 0 \\
0 & 0 & 1 & 0
\end{pmatrix},
\ 
S_-=
\begin{pmatrix}
0 & -3 & 0 & 0 \\
0 & 0 & -4 & 0 \\
0 & 0 & 0 & -3 \\
0 & 0 & 0 & 0
\end{pmatrix}.
\]
Let $\chi=\ul xX\ol\partial$, $\eta=\ul xH\ol\partial$, $\sigma_\pm=\ul xS_\pm\ol\partial$ be the corresponding vector fields.
Then 
\[
\begin{pmatrix}
\chi \\ \eta \\ \sigma_+ \\ \sigma_-
\end{pmatrix}
=A\cdot\ol\partial,
\quad
A=
\begin{pmatrix}
x_1 & x_2 & x_3 & x_4 \\
3x_1 & x_2 & -x_3 & -3x_4 \\
x_2 & x_3 & x_4 & 0 \\
0 & -3x_1 & -4x_2 & -3x_3
\end{pmatrix}.
\]
and $f=\frac{1}{2}\det(A)=3x_2^2x_3^2-6x_1x_3^3-8x_2^3x_4+18x_1x_2x_3x_4-9x_1^2x_4^2$ is irreducible.
By \cite[Lem.~1.9]{Sai80}, this implies that $f$ is free and $\chi,\eta,\sigma_+,\sigma_-$ is a basis of $\Der_f$.
By construction, $\DD_1=\C\times\mathfrak{sl}_2$ is not solvable.
So the statement of Proposition \ref{29} fails in dimension $n=4$.
\end{exa}

\section{LCT and Euler homogeneity}

Finally we return to the situation of the introduction.
Let $x\in D$ and choose a coordinate system $\ul x=(x_1,\dots,x_n)$ at $x$ defined in some Stein open neighbourhood $V_x$ of $x$.
Let $\V_x=(V_{x,i})_{i=1,\dots,n}$ be the Stein open covering of $V_x^*:=V_x\backslash\{x\}$ defined by $V_{x,i}:=\{x\in V_x\mid x_i\ne 0\}$ and $\U_x$ its restriction to $V_x\backslash D$.
The following explicit necessary condition for LCT for free divisors is due to F.J.~Castro--Jim\'enez et al.\ \cite{CNM96,CMNC02}.
Its proof relies on the comparison of the four spectral sequences arising from the two double complexes $\check\CC^q(\V_x,\Omega_{V_x^*}^p(\log D))$ and $\check\CC^q(\U_x,\Omega_{V_x\backslash D}^p)$.

\begin{thm}[(F.J.~Calder\'on--Moreno et al.)]\label{1}
Let $D\subseteq X$ be a free divisor and assume that LCT holds for $D$ in $V_x^*$ where $x\in X$.
Then LCT holds for $D$ at $x$ if and only if
\[
\xymat@C=15pt{
0\ar[r]&
\check\HH^{n-1}(V_x^*,\O_X)\ar[r]^-{\dd_1}&
\check\HH^{n-1}(V_x^*,\Omega_X^1(\log D))\ar[r]^-{\dd_2}&
\check\HH^{n-1}(V_x^*,\Omega_X^2(\log D))
}
\]
is an exact sequence.
\end{thm}

We shall combine Theorems \ref{6} and \ref{1} to prove Conjecture \ref{18} in dimension $n=3$, which is our main result Theorem \ref{30}.
We first give a more explicit description of the map $\dd_1$ in Theorem \ref{1}.

Since $\V_x$ is a Stein open covering of $V_x^*$, one can identify
\begin{align*}
H:&=\check\HH^{n-1}(V_x^*,\O_X)=\check\HH^{n-1}(\V_x,\O_X)\\
&=\C\bigl\{x_1,x_1^{-1},\dots,x_n,x_n^{-1}\bigr\}\Big/\sum_{i=1}^n\C\bigl\{x_1,x_1^{-1},\dots,x_i,\widehat{x_i^{-1}},\dots,x_n,x_n^{-1}\bigr\}
\end{align*}
where $\C\bigl\{x_1,x_1^{-1},\dots,x_n,x_n^{-1}\bigr\}$ denotes the ring of Laurent series in $x_1,\dots,x_n$.
Since the Stein open neighbourhoods of $x$ form a fundamental system of neighbourhoods of $x$, we may restrict our considerations to germs $X=(X,x)=(\C^n,0)$, $D=(D,x)=\VV(f)\subseteq X$ where $f\in\m_X$, and 
\[
\Der_X(-\log D)=\Der_f
\]
as in Definition \ref{25}.
In the following, we abbreviate $\m:=\m_X\subseteq\O_X=:\O$, $\Der:=\Der_X$, and $\Omega^\bullet:=\Omega_X^\bullet$.

Let $\ul\omega=(\omega_1,\dots,\omega_n)$ be a basis of $\Omega^1(\log D)$ and $\ul\delta=(\delta_1,\dots,\delta_n)$ its dual basis of $\Der(-\log D)$.
Via $\ul\omega:\O^n\cong\Omega^1(\log D)$ one can identify
\[
\xymat{
\dd_1=(\delta_1,\dots,\delta_n):H\ar[r]&H^n,\quad[g]\ar@{|->}[r]&\bigl([\delta_1(g)],\dots,[\delta_n(g)]\bigr).
}
\]
We shall only make use of the following consequence of Theorem \ref{1}.

\begin{cor}\label{27}
If LCT holds for a free divisor $D\subseteq X$ then $\ker\dd_1=0$.
\end{cor}

In the following, we abbreviate
\[
X':=(\C^{n-1},0),\quad X'':=(\C,0),\quad X:=X'\times X''
\]
and reduce the problem to the case where $D$ is not a product with a smooth factor.
We first note that freeness and strong Euler homogeneity are independent of smooth factors.

\begin{lem}\label{26}
Let $D\subseteq X$ be a divisor.
Then $D\cong D'\times X''$ for some divisor $D'\subseteq X'$ is equivalent to $\Der(-\log D)\not\subseteq\m\cdot\Der$.
In this case, $D$ is Euler homogeneous and each of the following properties is equivalent for $D$ and $D'$: strong Euler homogeneity, freeness, and Koszul freeness.
\end{lem}

\begin{proof}
This follows immediately from Lemma \ref{7}.
\end{proof}

By F.J.~Castro--Jim\'enez et al.\ \cite[Lem.~2.2.i,ii]{CNM96}, also LCT is independent of smooth factors.

\begin{lem}[(F.J.~Castro--Jim\'enez et al.)]\label{21}
Let $D'\subseteq X'$ be a divisor and $D=D'\times X''$.
Then LCT for $D'$ is equivalent to LCT for $D$.
\end{lem}

By Theorem \ref{31} and Lemmata \ref{26} and \ref{21}, we may assume from now on that $\II(D)=\langle f\rangle$ where $f$ is not a product as in Definition \ref{25}.
Then it suffices to prove that $f$ is Euler homogeneous if LCT holds for $D$ and the results in the preceding sections can be applied.

We shall use Corollary \ref{27} only in the following special case.

\begin{lem}\label{23}
Let $\delta\in\m\cdot\Der$ and $A\in\C^{n\times n}$ such that $\delta_0=\ul x A\ol\partial$.
Then 
\[
\delta\Bigl[\frac{1}{x_1\cdots x_n}\Bigr]=\Bigl[\frac{\tr A}{x_1\cdots x_n}\Bigr]\in H.
\]
\end{lem}

\begin{proof}
This follows immediately from the definition of $H$.
\end{proof}

From now on, let $n=3$ and abbreviate
\[
x,y,z:=x_1,x_2,x_3\quad\text{ and }\quad\partial_x,\partial_y,\partial_z:=\partial_{x_1},\partial_{x_3},\partial_{x_3}.
\]
We assume that $f$ is not Euler homogeneous and claim that LCT does not hold for $D$.
By Corollary \ref{12}, there is an $\wh\O$--basis $\sigma_1,\dots,\sigma_s,\nu_1,\dots,\nu_{n-s}$ of $\Der_f\wh\O$ as in Theorem \ref{6} and, by Corollary \ref{22}, we may assume that $\sigma_i(f)=0$ for all $i=1,\dots,s$.
There are the following cases:

\begin{asparaenum}[{\bf Case} I:]

\item $s=0$.
Then the claim follows from Corollary \ref{27} and Lemma \ref{23} using a truncated coordinate change in Theorem \ref{6}.

\item $s=1$ and $\sigma=\sigma_1=ax\partial_x+by\partial_y+cz\partial_z$.

\begin{asparaenum}[{\em Subcase} a:]

\item $a\ne0$ and $b,c=0.$
In this situation, $f$ is annihilated by $\sigma = a x\partial_x$ and hence $\partial_x\in\Der_f$ in contradiction to our assumption that $f$ is not a product.

\item $a,b\ne0$ and $c=0$.
Then $\sigma(f)=0$ implies $ab<0$ and $f=\sum_{ia+jb=0}a_{i,j}(z)x^iy^j$.
Since $f\in\langle x,y\rangle$ by Definition \ref{15}, $a_{0,0}(z)=0$. 
But $f$ being reduced implies $a_{1,1}(z)\neq 0$ which forces $a=-b$.
Then the claim follows from Corollary \ref{27} and Lemma \ref{23} using a truncated coordinate change in Theorem \ref{6}.

\item $a,b,c\ne0$.
A truncated coordinate change in Theorem \ref{6} yields the existence of a (convergent)
 $\delta\in\Der_f\O$ such that $\delta_0=\sigma$.
Then $\delta$ vanishes only at the origin and is tangent to the $1$--dimensional smooth part of $\Sing D$.
This implies that the logarithmic characteristic subvariety $\LL_X(-\log D)\subseteq\TT_X^*$ \cite[Def.~3.15]{Sai71} has minimal dimension $n=3$.
(In the language of K.~Saito, the existence of $\delta$ above implies
that the logarithmic stratification of $X$ consists only 
of holonomic strata \cite[Def.~3.3, 3.8]{Sai71} and hence $\LL_X(-\log D)$ has only holonomic components \cite[Def.~3.17, Prop.~3.18]{Sai71}.)
But then $D$ is Koszul free by \cite[Cor.~1.9]{CN02} and LCT does not hold
for $D$ by Corollary \ref{32}.

\end{asparaenum}

\item $s\geq 2$.
In this situation, there are two linearly independent $\sigma_1,\sigma_2\in\Der_f\wh\O$ and hence the Newton diagram of $f$ is contained in a $1$--dimensional vector space.
Then there is a monomial $x^iy^jz^k$ such that $f=u\cdot x^iy^jz^k$ for some $u\in\wh\O^*$.
But this means that $D$ is Euler homogeneous in contradiction to our assumption. 

\end{asparaenum}

Finally we have proved our main result Theorem \ref{30}.
There is also a simple proof of Theorem \ref{31} using Theorem \ref{6}, Corollary \ref{27}, and Lemma \ref{23} as above.

\bibliographystyle{amsalpha}
\bibliography{fslvf}

\providecommand{\bysame}{\leavevmode\hbox to3em{\hrulefill}\thinspace}
\providecommand{\MR}{\relax\ifhmode\unskip\space\fi MR }
\providecommand{\MRhref}[2]{%
  \href{http://www.ams.org/mathscinet-getitem?mr=#1}{#2}
}
\providecommand{\href}[2]{#2}
\begin{thebibliography}{CMNC02}

\bibitem[AA88]{AA88}
D.V. Anosov and I.V. Arnol'd (eds.), \emph{Dynamical systems {I}}, EMS, vol.~1,
  Springer, 1988.

\bibitem[CMNC02]{CMNC02}
F.J. {Calder{\'o}n--Moreno}, D.~Mond, L.~{Narv{\'a}ez--Macarro}, and F.J.
  {Castro--Jim{\'e}nez}, \emph{Logarithmic cohomology of the complement of a
  plane curve}, Comm. Math. Helv. \textbf{77} (2002), 24--38.

\bibitem[CN02]{CN02}
F.J. {Calder{\'o}n--Moreno} and L.~{Narv{\'a}ez--Macarro}, \emph{The module
  $\mathcal{D}f^s$ for locally quasi--homogeneous free divisors}, Comp. Math.
  \textbf{134} (2002), 59--74.

\bibitem[CN05]{CN05}
\bysame, \emph{Dualit{\'e} et comparaison sur les complexes de {de Rham}
  logarithmiques par rapport aux diviseurs libres}, Ann. Inst. Fourier Grenoble
  \textbf{55} (2005), no.~1, 47--75.

\bibitem[CNM96]{CNM96}
F.J. {Castro--Jim{\'e}nez}, L.~{Narv{\'a}ez--Macarro}, and D.~Mond,
  \emph{Cohomology of the complement of a free divisor}, Trans. Amer. Math.
  Soc. \textbf{348} (1996), no.~8, 3037--3049.

\bibitem[CU05]{CU05}
F.J. {Castro--Jim\'enez} and J.M. {Ucha--Enr\'iquez}, \emph{{Logarithmic
  Comparison Theorem} and some {Euler} homogeneous free divisors}, Proc. Amer.
  Math. Soc. \textbf{133} (2005), 1417--1422.

\bibitem[Del71]{Del71}
P.~Deligne, \emph{Th\'eorie de {Hodge}, {II}}, Publ. Math. I.H.E.S \textbf{40}
  (1971), 5--57.

\bibitem[Gro66]{Gro66}
A.~Grothendieck, \emph{On the de {Rham} cohomology of algebraic varieties},
  Publ. Math. de l'I.H.E.S. \textbf{29} (1966), 95--103.

\bibitem[Jac62]{Jac62}
N.~Jacobson, \emph{Lie algebras}, Interscience, 1962.

\bibitem[Sai71]{Sai71}
K.~Saito, \emph{{Quasihomogene isolierte Singularit\"aten von Hyperfl\"achen}},
  Inv. Math. \textbf{14} (1971), 123--142.

\bibitem[Sai80]{Sai80}
\bysame, \emph{Theory of logarithmic differential forms and logarithmic vector
  fields}, J. Fac. Sci. Univ. Tokyo \textbf{27} (1980), 265--291.

\bibitem[Tor04]{Tor04}
T.~Torrelli, \emph{On meromorphic functions defined by a differential system of
  order 1}, Bull. Soc. Math. France \textbf{132} (2004), 591--612.

\end{thebibliography}

\end{document}